\documentclass[12pt]{amsart}

\newtheorem{theorem}{Theorem}[section]
\newtheorem{lemma}[theorem]{Lemma}
\newtheorem{proposition}[theorem]{Proposition}
\newtheorem{corollary}[theorem]{Corollary}
\newtheorem{definition}[theorem]{Definition}
\newtheorem{remark}[theorem]{Remark}

\numberwithin{equation}{section}

\begin{document}

\baselineskip=15.5pt

\title{On the stratified vector bundles}

\author[I. Biswas]{Indranil Biswas}

\address{School of Mathematics, Tata Institute of Fundamental
Research, Homi Bhabha Road, Bombay 400005, India}

\email{indranil@math.tifr.res.in}



\date{}

\begin{abstract}

The stratified vector bundles on a smooth variety defined
over an algebraically closed field $k$ form a neutral Tannakian
category over $k$. We investigate the affine group--scheme
corresponding to this neutral Tannakian category.

\end{abstract}

\maketitle

\section{Introduction}

Let $X$ be an irreducible smooth variety defined over an
algebraically closed field $k$ of positive characteristic.
Gieseker in \cite{Gi} introduced the notion of a stratified
vector bundle over $X$. These objects are analogs of the complex
algebraic vector bundles equipped with an integrable connection.
More precisely, a stratified vector bundle over $X$ is a pair
$(E\, ,\nabla)$, where $E$ is an algebraic vector bundle
over $X$, and 
$$
\nabla\, :\, {\mathcal D}(X)\, \longrightarrow\,
{\mathcal E}nd_k(E)
$$
is an ${\mathcal O}_X$--linear ring homomorphism, where
${\mathcal D}(X)$ is the sheaf of differential operators sending
${\mathcal O}_X$ to itself, and ${\mathcal E}nd_k(E)$ is the
sheaf of $k$--linear endomorphisms of $E$.

The stratified vector bundles over $X$
form a rigid abelian $k$--linear tensor category;
this category will be denoted by ${\mathcal C}(X)$.
After we fix a $k$--rational
point $x_0$ of $X$, this category ${\mathcal C}(X)$ gets the
natural fiber functor that sends any 
stratified vector bundle $(E\, ,\nabla)$ to the fiber of $E$
over the point $x_0$. Therefore, ${\mathcal C}(X)$
equipped with this fiber functor is a neutral Tannakian
category over $k$.
Consequently, we get an affine group scheme over $k$, which
will be called the stratified group--scheme of $X$. This
stratified group--scheme depends only on the pointed
variety $(X\, ,x_0)$, and it is denoted by
${\mathcal S}(X,x_0)$.

The group--scheme
${\mathcal S}(X,x_0)$ has a surjective homomorphism
to the \'etale fundamental group $\pi^{et}(X,x_0)$
(Proposition \ref{prop1}). Hence ${\mathcal S}(X,x_0)$
in general is not of finite type. The Frobenius
morphism of the group--scheme ${\mathcal S}(X,x_0)$
turns out to be an automorphism (Corollary \ref{cor.fi}).
We also show that ${\mathcal S}(X,x_0)$ is different from
the fundamental group--scheme constructed by Nori
(the details are in Section \ref{ex.}).

If $X$ and $Y$ are irreducible smooth proper varieties
over $k$, then we show that
$$
{\mathcal S}(X\times Y,(x_0,y_0))\, =\,
{\mathcal S}(X,x_0)\times {\mathcal S}(Y,y_0)
$$
(see Theorem \ref{thm1}).

For any nonnegative integer $n$, let $X_n$ denote the
base change of $X$ using the homomorphism $k\, \longrightarrow
\, k$ defined by $\lambda\, \longmapsto\, \lambda^{p^n}$,
where $p$ is the characteristic of $k$.
In \cite{Gi}, Gieseker defined a flat vector bundle over $X$
to be sequence
of pairs $\{E_n\, ,\sigma_n\}_{n\geq 0}$, where $E_n$ is a
vector bundle over $X_n$, and
$$
\sigma_n\, :\, F^*_X E_{n+1}\, \longrightarrow\, E_{n}
$$
is an isomorphism of vector bundles, where $F_X$ is the
relative Frobenius morphism. As shown in \cite{Gi},
there is a canonical equivalence of categories between the
category of flat vector bundles on $X$ and the category
${\mathcal C}(X)$ of stratified vector bundles.

Let $\mathcal G$ be an affine group--scheme
defined over the field $k$.
A flat principal $\mathcal G$--bundle over $X$ is a
sequence of pairs $\{E_n\, ,\sigma_n\}_{n\geq 0}$, where $E_n$
is a principal $\mathcal G$--bundle over $X_n$ and
$$
\sigma_n\, :\, F^*_X E_{n+1}\, \longrightarrow\, E_{n}
$$
is an isomorphism of principal $\mathcal G$--bundles.
A flat principal $\mathcal G$--bundle $\{E_n\, ,
\sigma_n\}_{n\geq 0}$ is also called a flat structure
on the principal $\mathcal G$--bundle $E_0$.

We show that there is a tautological principal
${\mathcal S}(X,x_0)$--bundle $E_{{\mathcal S}(X,x_0)}$ on
$X$. This principal ${\mathcal S}(X,x_0)$--bundle
$E_{{\mathcal S}(X,x_0)}$ is equipped with a
tautological flat structure. (See Section \ref{se5.1}
and Lemma \ref{lem4}.) We will describe below a universal
property of $E_{{\mathcal S}(X,x_0)}$.

Let $\{E_n\, ,\sigma_n\}_{n\geq 0}$ be a flat principal
$\mathcal G$--bundle over $X$, where
$\mathcal G$ is any affine group--scheme over $k$. Let
${\rm Ad}(E_0)$ be the adjoint bundle of the principal
$\mathcal G$--bundle $E_0$; so ${\rm Ad}(E_0)$ is a group--scheme
over $X$. Fixing a $k$--rational point $\widehat{x}_0$
in the fiber $(E_0)_{x_0}$ over $x_0$, we get an isomorphism
of group--schemes
$$
\rho_0\, :\, \text{Ad}(E_0)_{x_0}\, \longrightarrow
\, {\mathcal G}\, .
$$

The tautological flat bundle $E_{{\mathcal S}(X,x_0)}$
over $X$ has the following universal property (see
Proposition \ref{prop2}):

There is homomorphism of group schemes
$$
\rho\, :\, {\mathcal S}(X,x_0)\, \longrightarrow\,
{\rm Ad}(E_0)_{x_0}
$$
canonically associated to the flat principal
$\mathcal G$--bundle $\{E_n\, ,\sigma_n\}_{n\geq 0}$. This
homomorphism $\rho$ has the property that the
flat principal $\mathcal G$--bundle $\{E_n\, ,
\sigma_n\}_{n\geq 0}$ coincides with the one obtained by
extending the structure group of the flat principal ${\mathcal 
S}(X,x_0)$--bundle $E_{{\mathcal S}(X,x_0)}$
using the homomorphism $\rho_0\circ \rho$, where $\rho_0$
is the above isomorphism.

\section{Stratified vector
bundles and flat vector bundles}\label{sec2}

\subsection{Preliminaries}\label{sec2.1}

Let $k$ be an algebraically closed field of characteristic
$p$, with $p\, >\, 0$.
Let $X$ be an irreducible
smooth variety defined over $k$. For any
integer $n\, \geq\, 0$, let
$$
{\mathcal D}^n(X)\, :=\, \text{Diff}^n({\mathcal O}_X\, ,
{\mathcal O}_X)
$$
be the sheaf of differential operators of order at most $n$
mapping ${\mathcal O}_X$ to itself. The direct limit
\begin{equation}\label{DX}
{\mathcal D}(X)\,:=\,
\lim_{\stackrel{\longrightarrow}{n}}{\mathcal D}^n(X)
\end{equation}
is called the \textit{sheaf of differential operators} on $X$.

We recall from \cite{Gi} the definition of a
stratified sheaf (see \cite[page 2, Definition 0.3]{Gi}).

\begin{definition}\label{def1}
{\rm A} stratified sheaf {\rm on $X$ is an
${\mathcal O}_X$--coherent sheaf $E$
on $X$ together with a ring homomorphism}
$$
\nabla\, :\, {\mathcal D}(X)\, \longrightarrow\, {\mathcal E}nd_k
(E)
$$
{\rm which is ${\mathcal O}_X$--linear.}

\begin{itemize}
\item {\rm A stratified sheaf $(E\, ,\nabla)$ with the
${\mathcal O}_X$--module $E$ locally
free is called a} stratified vector bundle.

\item {\rm For a stratified sheaf $(E\, ,\nabla)$, the
homomorphism
$\nabla$ is called a} stratification {\rm on $E$.}
\end{itemize}
\end{definition}

\begin{lemma}\label{lem1}
For any stratified sheaf on $X$, the underlying
${\mathcal O}_X$--coherent sheaf is locally free.
\end{lemma}

See \cite[\S~2, page 21, Proposition 2.16]{BO} for a proof
of the above lemma.

A \textit{homomorphism} from a stratified vector bundle
$(E\, , \nabla)$ to a stratified vector bundle $(E'\, ,
\nabla')$ is a homomorphism of ${\mathcal O}_X$--coherent 
sheaves
$$
f\, :\, E\, \longrightarrow\, E'
$$
that intertwines the actions
of ${\mathcal D}(X)$ on $E$ and $E'$, or in other words,
$f(\nabla(D)(s))\, =\, \nabla'(D)(f(s))$, where
$D$ is a locally defined section of ${\mathcal D}(X)$, and
$s$ is a locally defined section of $E$. (In \cite{Gi}, a
homomorphism from $(E\, , \nabla)$ to $(E'\, , \nabla')$ is also
called a \textit{horizontal homomorphism} from $E$ to $E'$.)

Using Lemma \ref{lem1} it is easy to check that the
stratified vector bundles on $X$ form an abelian category.

For any nonnegative integer $n$, let
$X_n\,=\, X\times_k k$ be the base change of
the variety $X$ using the homomorphism $k\,
\longrightarrow\, k$ defined by $t\, \longmapsto\, t^{p^n}$.
So in particular $X\, =\, X_0$.
Since the field $k$ is algebraically closed, in particular it
is perfect, $X_n$ is isomorphic to $X$ for all $n$.
For any $n\, \geq\, 0$, let
$$
F_X\, :\, X_n\, \longrightarrow\, X_{n+1}
$$
be the relative Frobenius morphism. We should clarify
that this is an abuse of notation since the domain of $F_X$
is not fixed and it depends on $n$. However from the context
the domain will generally be clear, and whenever there is any
ambiguity, we will clearly specify the domain. For any positive
integer $d$, let
$$
F^d_X \, :=\, \overbrace{F_X\circ\cdots\circ F_X}^{d\mbox{-}\rm{times}}
\, :\, X_n\, \longrightarrow\, X_{n+d}
$$
be the $d$--fold iteration of $F_X$. For notational convenience,
$F^0_X$ will denote the identity map of $X_n$.
(See \cite[Section 4.1]{Ra} for details.)

We recall from \cite{Gi} the definition of a flat vector bundle
(see \cite[page 3, Definition 1.1]{Gi}).

\begin{definition}\label{def2}
{\rm A} flat {\rm sheaf over $X$ is a sequence
of pairs $\{E_n\, ,\sigma_n\}_{n\geq 0}$, where $E_n$ is an
${\mathcal O}_X$--coherent sheaf on $X_n$ and}
$$
\sigma_n\, :\, F^*_X E_{n+1}\, \longrightarrow\, E_{n}
$$
{\rm is an isomorphism of ${\mathcal O}_X$--coherent sheaves.}

\begin{itemize}
\item {\rm A flat sheaf $\{E_n\, ,\sigma_n\}_{n\geq 0}$ with all
${\mathcal O}_X$--modules $E_n$ locally free is called a}
flat vector bundle.

\item {\rm A flat sheaf $\{E_n\, ,\sigma_n\}_{n\geq 0}$
will also be called a} flat structure on $E_0$.
\end{itemize}
\end{definition} 

A \textit{homomorphism} from a flat sheaf
$\{E_n\, ,\sigma_n\}_{n\geq 0}$ to a flat sheaf
$\{E'_n\, ,\sigma'_n\}_{n\geq 0}$ is a homomorphism
of ${\mathcal O}_X$--coherent sheaves
$$
\tau_n\, : \, E_n\, \longrightarrow\, E'_n
$$
for each $n\, \geq\, 0$ such that the two homomorphisms $\sigma'_n
\circ (F^*_X \tau_{n+1})$ and $\tau_n\circ\sigma_n$ from
$F^*_X E_{n+1}$ to $E'_n$ coincide. (In \cite{Gi}, a
homomorphism from $\{E_n\, ,\sigma_n\}_{n\geq 0}$ to
$\{E'_n\, ,\sigma'_n\}_{n\geq 0}$ is called a \textit{horizontal
homomorphism} from $E_0$ to $E'_0$.)

A theorem due to Katz identifies stratified vector bundles
with flat vector bundles (see \cite[page 4, Theorem 1.3]{Gi}).
Using Lemma \ref{lem1}, from the proof of Theorem 1.3 in \cite{Gi}
it follows immediately that if $\{E_n\, ,\sigma_n\}_{n\geq 0}$
is a flat sheaf, then the ${\mathcal O}_X$--coherent sheaf $E_n$ is 
locally free for each $n$. Combining
this together with Lemma \ref{lem1} we conclude that
the flat vector bundles over $X$ form an abelian category.

\subsection{An affine group--scheme}

The tensor product of two stratified vector bundles
$(E\, , \nabla)$ and $(E'\, , \nabla')$ is constructed
as follows. Take any locally defined vector field
$\xi$ on $X$, and let $u$ and $v$ be locally defined
sections of $E$ and $E'$ respectively. Define
$$
\widehat{\nabla}(\xi)(u\otimes v) \, :=\, (\nabla
(\xi)(v))\otimes w + v \otimes (\nabla' (\xi)(w))\, ,
$$
which is a locally defined section of $E\bigotimes E'$.
This operator $\widehat{\nabla}$ extends to an
${\mathcal O}_X$--linear ring homomorphism
$$
\widehat{\nabla}\, :\, {\mathcal D}(X)\, \longrightarrow\,
{\mathcal E}nd_k (E\otimes E')\, .
$$
The tensor product $(E\, , \nabla)\bigotimes (E'\, , \nabla')$
is defined to be $(E\bigotimes E'\, , \widehat{\nabla})$.

We note that if $\{E_n\, ,\sigma_n\}_{n\geq 0}$ and
$\{E'_n\, ,\sigma'_n\}_{n\geq 0}$ are the flat vector
bundle over $X$ corresponding to the stratified vector bundles
$(E\, , \nabla)$ and $(E'\, , \nabla')$ respectively, then
the flat vector bundle corresponding to the tensor product $(E\, , 
\nabla)\bigotimes (E'\, , \nabla')$ is
$$
\{E_n\otimes E'_n\, ,\sigma_n\otimes \sigma'_n\}_{n\geq 0}\, .
$$

Let ${\mathcal C}(X)$ denote the category of
stratified vector bundles on $X$. We already noted
that ${\mathcal C}(X)$ is an abelian category.
In fact, ${\mathcal C}(X)$ is a rigid abelian $k$--linear
tensor category (see \cite[page 112, Definition 1.7]{DM} for
the definition of a rigid abelian $k$--linear tensor category).

Let $\text{Vect}(k)$
denote the category of finite dimensional $k$--vector
spaces. Fix a $k$--rational point $x_0\,\in\, X$.

We have a fiber functor
\begin{equation}\label{e1}
T_{x_0}\, :\, {\mathcal C}(X)\, \longrightarrow\,
\text{Vect}(k)
\end{equation}
that sends a stratified vector bundle $(E\, , \nabla)$ over $X$
to its fiber $E_{x_0}$ over the base point $x_0$.
Equipped with this
fiber functor $T_{x_0}$, the category ${\mathcal C}(X)$
of stratified vector bundles becomes a neutral Tannakian
category over $k$ (see \cite[page 138, Definition 2.19]{DM}).

Therefore, the neutral Tannakian category
$({\mathcal C}(X)\, , T_{x_0})$, where
$T_{x_0}$ is defined in Eq. \eqref{e1}, produces an affine
group--scheme over $k$ \cite[page 130, Theorem
2.11]{DM}, \cite[Theorem 1.1]{No}, \cite[Theorem 1]{Sa}.

\begin{definition}\label{def3}
{\rm The group scheme over $k$ given by the neutral Tannakian
category $({\mathcal C}(X)\, , T_{x_0})$ over $k$
will be called the} stratified group--scheme {\rm of $X$.}

{\rm The stratified group--scheme of $X$ will be denoted by
${\mathcal S}(X,x_0)$.}
\end{definition}

Let $Y$ be an irreducible
smooth variety defined over $k$ and
$$
\varphi\, :\, X\, \longrightarrow\, Y
$$
a morphism. Set $y_0\, :=\, \varphi(x_0)\, \in\, Y$.
If $(E\, ,\nabla)$ is a stratified vector bundle over $Y$,
then the pull back $(\varphi^*E\, ,\varphi^*\nabla)$ is a
stratified vector bundle over $X$. It is easy to that if 
$$
f\, :\, (E\, ,\nabla)\, \longrightarrow\, (E'\, ,\nabla')
$$
is a homomorphism between stratified vector bundles over $Y$,
then
$$
\varphi^*f\, :\, (\varphi^* E\, ,\varphi^*\nabla)\,
\longrightarrow\, (\varphi^*E'\, ,\varphi^*\nabla')
$$
is also a homomorphism between stratified vector bundles.
Consequently, the morphism $\varphi$ from $X$ to $Y$
induces a homomorphism of group--schemes
\begin{equation}\label{e-v}
\varphi^*\, :\, {\mathcal S}(X,x_0)\,\longrightarrow\,
{\mathcal S}(Y,y_0)\, .
\end{equation}

Assume that $X$ is proper. Let $(E\, ,\nabla)$ be
a stratified vector bundle over $X\times Z$, where
$Z$ is an irreducible smooth variety defined over $k$.
Let
$$
\psi\, :\, X\times Z\, \longrightarrow\, Z
$$
be the natural projection.
Since $X$ is proper, the direct image $R^i\psi_{*}E$ is
a coherent sheaf on $Z$ for each $i\, \geq\, 0$
\cite[page 116, Th\'eor\`eme 3.2.1]{Gr}. The
pull back $\psi^*{\mathcal D}(Z)$ is canonically a
subsheaf of ${\mathcal D}(X\times Z)$ (see
Eq. \eqref{DX} for definition). Using this inclusion
of $\psi^*{\mathcal D}(Z)$ in ${\mathcal D}(X\times Z)$,
the homomorphism
$$
\nabla\, :\, {\mathcal D}(X\times Z)\, \longrightarrow\,
{\mathcal E}nd_k(E)
$$
produces a homomorphism
$$
\nabla'\, :\, {\mathcal D}(Z)\, \longrightarrow\,
{\mathcal E}nd_k(R^i\psi_{*}E)
$$
for each $i\, \geq\, 0$. It is straight--forward to check that
\begin{equation}\label{d.c}
(R^i\psi_{*}E\, , \nabla')
\end{equation}
is a stratified sheaf on $Z$. From Lemma \ref{lem1}
it now follows that the ${\mathcal O}_Z$--coherent
sheaf $R^i\psi_{*}E$ is locally free.

\begin{remark}\label{rem1}
{\rm Let $X$ be an irreducible smooth variety over $k$ and
$U\, \subset\, X$ a nonempty Zariski open subset such that
the complement $X\setminus U$ is of codimension at least
two. Fix a $k$--rational point $x_0$ of $U$. The inclusion
map of $(U\, ,x_0)$ in $(X\, ,x_0)$ induces an isomorphism}
$$
{\mathcal S}(U,x_0)\, \stackrel{\sim}{\longrightarrow}\,
{\mathcal S}(X,x_0)
$$
{\rm of stratified group--schemes. This follows
immediately from \cite[page 21, Theorem 3.14]{Gi}.}
\end{remark}

\section{Some properties of the group--scheme}

\subsection{The \'etale fundamental group as a quotient group}

Throughout this section we will assume $X$ to be an
irreducible smooth variety proper over $k$.

A vector bundle $E$ over $X$ is called \textit{\'etale
trivializable} if there is an algebraic \'etale Galois covering
$$
\gamma\, :\, Y\, \longrightarrow\, X
$$
such that the pull back $\gamma^*E$ is trivializable. Therefore,
a vector bundle $E$ of rank $n$ over $X$ is \'etale trivializable
if and only if there is a representation
$$
\rho\, :\, \pi^{et}(X,x_0)\, \longrightarrow\,
\text{GL}(n,k)\, ,
$$
where $\pi^{et}(X,x_0)$ is the \'etale fundamental group of $X$
with base point $x_0$, such that the associated vector bundle
$V_\rho$ over $X$ is isomorphic to $E$.

We briefly recall a Tannakian description of $\pi^{et}(X,x_0)$.

The \'etale trivializable vector bundles over $X$ from a
rigid abelian $k$--linear
tensor category ${\mathcal E}(X)$. For any
$E\, ,E'\, \in\, {\mathcal E}(X)$,
the homomorphisms from $E$ to $E'$ are defined to be all
${\mathcal O}_X$--linear homomorphisms from the
vector bundle $E$ to $E'$. The direct sum, tensor
product and duals are defined in the obvious way.
We have the fiber functor
\begin{equation}\label{Tu0}
T^0_{x_0}\, :\, {\mathcal E}(X)\, \longrightarrow\,
\text{Vect}(k)
\end{equation}
that sends any \'etale trivializable vector bundle
$E$ to its fiber $E_{x_0}$ over the point $x_0$.
The group--scheme given by the
neutral Tannakian category $({\mathcal E}(X)\, ,
T^0_{x_0})$ over $k$
is identified with the \'etale fundamental group
$\pi^{et}(X,x_0)$.

Take any representation
$$
\rho\, :\, \pi^{et}(X,x_0)\, \longrightarrow\,
\text{GL}(n,k)\, .
$$
Let $V_\rho$ be the vector bundle of rank $n$ over $X$
associated to $\rho$. This $V_\rho$ is equipped with a
canonical stratification (see \cite[page 7]{Gi} for the details).

\begin{lemma}\label{lem2}
Let $E$ and $F$ be two \'etale trivializable vector bundles
over $X$. Then any ${\mathcal O}_X$--linear homomorphism from
the vector bundle $E$ to $F$ intertwines the
${\mathcal D}(X)$--module structures defining the
stratifications on $E$ and $F$.
\end{lemma}

\begin{proof}
Let $\gamma\, :\, Y\, \longrightarrow\, X$ be an algebraic \'etale
Galois covering with $Y$ connected such that $\gamma^*E$ is
trivializable. Similarly, take an algebraic \'etale Galois covering
$$
\delta\, :\, Z \, \longrightarrow\, X
$$
with $Z$ connected such that $\delta^*F$ is trivializable.
Fix a connected component
$$
M\, \subset\, Y\times_X Z
$$
of the fiber product. Let
\begin{equation}\label{phi}
\phi\, :\, M \, \longrightarrow\, X
\end{equation}
be the natural projection obtained by restricting the
map $\gamma\times \delta$. Since $\gamma^*E$ and
$\delta^*F$ are trivializable, it follows immediately that
$\phi^*E$ and $\phi^*F$ are also trivializable.

Let
$$
h\,:\, E\, \longrightarrow\, F
$$
be any ${\mathcal O}_X$--linear homomorphism. Let
\begin{equation}\label{phi2}
\phi^*h\, :\, \phi^*E\, \longrightarrow\, \phi^*F
\end{equation}
be the pull back of $h$ to $M$, where $\phi$ is the
projection in Eq. \eqref{phi}. We noted earlier that
both $\phi^*E$ and $\phi^*F$ are trivializable. Fix
trivializations of $\phi^*E$ and $\phi^*F$. With
respect to these trivializations, the homomorphism
$\phi^*h$ in Eq. \eqref{phi2} is given by a morphism from
$M$ to the variety of $n\times m$ matrices
\begin{equation}\label{H}
H\, :\, M\, \longrightarrow\, \text{M}_{n\times m}(\mathbb C)\, ,
\end{equation}
where $m\, =\, \text{rank}(E)$ and $n\, =\, \text{rank}(F)$.

The variety $M$ is proper over $k$ because $X$ is proper.
Hence the morphism $H$ in Eq. \eqref{H} must be a constant
one. This immediately implies that $h$ intertwines the
${\mathcal D}(X)$--module structures of $E$ and $F$. This
completes the proof of the lemma.
\end{proof}

In view of Lemma \ref{lem2}, we obtain a
functor between the neutral Tannakian categories over $k$
\begin{equation}\label{cF}
{\mathcal F}\, :\, ({\mathcal E}(X)\, ,T^0_{x_0})\,
\longrightarrow\, ({\mathcal C}(X)\, , T_{x_0})\, ,
\end{equation}
where $T_{x_0}$ and $T^0_{x_0}$ are constructed in
Eq. \eqref{e1} and Eq. \eqref{Tu0} respectively.
This functor ${\mathcal F}$
produces a homomorphism of group--schemes over $k$
\begin{equation}\label{h.e.f.}
\theta\,:\,{\mathcal S}(X,x_0)\,\longrightarrow\,\pi^{et}(X,x_0)\, ,
\end{equation}
where ${\mathcal S}(X,x_0)$ is constructed in Definition \ref{def3}.

If $\rho'\, :\, \pi^{et}(X,x_0)\, \longrightarrow\,
\text{GL}(n,k)$ is another representation which
is not isomorphic to $\rho$, then the stratified vector
bundle $V_{\rho'}$ over $X$ associated to $\rho'$ is not
isomorphic to the stratified
vector bundle $V_{\rho}$ associated to $\rho$
\cite[page 7, Proposition 1.9]{Gi}. From this it follows
that the image $\theta({\mathcal S}(X,x_0))$ of the homomorphism
$\theta$ in Eq. \eqref{h.e.f.} is not contained in any
normal proper subgroup of $\pi^{et}(X,x_0)$.
In fact, a stronger statement holds. The homomorphism
$\theta$ is surjective, as shown by the following proposition.

\begin{proposition}\label{prop1}
The homomorphism $\theta$ in Eq. \eqref{h.e.f.} is faithfully
flat. In particular, $\theta$ is surjective.
\end{proposition}

\begin{proof}
We will use the criterion in \cite[page 139, Proposition 2.21(a)]{DM}
for a homomorphism of group--schemes to be faithfully flat.
The functor $\mathcal F$ in Eq. \eqref{cF} is
evidently fully faithful. Therefore, to complete the proof
using the criterion in \cite[Proposition 2.21(a)]{DM} we need
to show the following.

Let $E$ be an \'etale trivializable
vector bundle over $X$. Let
$$
\nabla\, :\, {\mathcal D}(X)\, \longrightarrow\, {\mathcal E}nd_k
(E)
$$
be the canonical stratification on $E$. Take any coherent subsheaf
$$
F\, \subset\, E
$$
preserved by $\nabla$, i.e.,
\begin{equation}\label{e3}
\nabla ({\mathcal D}(X))(F)\, \subset\, F\, .
\end{equation}
Then $F$ is also \'etale trivializable.

To prove that $F$ is \'etale trivializable, take an algebraic
\'etale Galois covering
$$
\gamma\, :\, Y\, \longrightarrow\, X
$$
with $Y$ connected such that $\gamma^*E$ is
trivializable. The condition in Eq. \eqref{e3} ensures
that the ${\mathcal O}_Y$--coherent subsheaf
$$
\gamma^*F\, \subset\, \gamma^*E
$$
is preserved by the natural connection on the trivializable
vector bundle $\gamma^*E$. This implies that
\begin{itemize}
\item $\gamma^*F$ is subbundle of $\gamma^*E$, and
\item the morphism $Y\, \longrightarrow\, \text{Gr}(r,\gamma^*E)$
associated to the subbundle $\gamma^*F$ is a constant one, where
$r\, =\, \text{rank}(F)$, and $\text{Gr}(r,\gamma^*E)$ is the
Grassmann bundle over $Y$ parametrizing $r$ dimensional subspaces
in the fibers of $\gamma^*E$.
\end{itemize}
Therefore, the subbundle $\gamma^*F$ is trivializable. In
particular, $F$ is \'etale trivializable. This completes
the proof of the proposition.
\end{proof}

\begin{remark}\label{rem2}
{\rm There are examples of stratified vector bundles on
smooth projective curves which are not \'etale trivializable;
see \cite[page 100]{Gi0}. Therefore, the homomorphism
$\theta$ in Eq. \eqref{h.e.f.} is not a closed embedding
in general.}
\end{remark}

\subsection{An essentially finite vector bundle which is
not flat}\label{ex.}

In \cite{No0}, Nori introduced the notion of
an essentially finite
vector bundle over a projective variety $X$. He showed that the
essentially finite vector bundles
form a neutral Tannakian category once a point $x_0$
of $X$ is fixed. The group--scheme $\pi(X,x_0)$
associated to this neutral Tannakian category given by the
essentially finite vector bundles is known as the
\textit{fundamental group--scheme}.

We will give an example of an essentially finite vector bundle
that is not flat. Such an example shows that the stratified
group--scheme ${\mathcal S}(X,x_0)$ does not coincide with the 
fundamental group--scheme $\pi(X,x_0)$ in general.

Let $X$ be a supersingular elliptic curve defined over an
algebraically closed field of positive characteristic.
This means that the pull back homomorphism
$$
F^*_X \,:\, H^1(X,\, {\mathcal O}_X)\, \longrightarrow\,
H^1(X,\, {\mathcal O}_X)
$$
vanishes. Fix a nonzero element
\begin{equation}\label{om.}
\omega\, \in\, H^1(X,\, {\mathcal O}_X)\, .
\end{equation}
Let
\begin{equation}\label{sh.ex.}
0\, \longrightarrow\, {\mathcal O}_X \, \longrightarrow\,
E \, \stackrel{f}{\longrightarrow}\,
{\mathcal O}_X \, \longrightarrow\, 0
\end{equation}
be the extension given by $\omega$.

The cohomology class $F^*_X\omega \, \in\,
H^1(X,\, {\mathcal O}_X)$ vanishes because $X$ is
supersingular. Hence the short exact sequence in
Eq. \eqref{sh.ex.} splits. Therefore,
$F^*_XE \, =\, {\mathcal O}_X\bigoplus {\mathcal O}_X$.
This implies that the vector bundle $E$ is essentially
finite (see \cite[pp. 552--553, Proposition 2.3]{BH}).

We will show that there is no vector bundle $V$ over $X$
such that $F^*_X V\, =\, E$. To prove this by contradiction,
let $V$ be a vector bundle over $X$ such that
\begin{equation}\label{sh.ex2.}
F^*_X V\, =\, E\, .
\end{equation}
Since $E$ is indecomposable of degree zero, from Eq.
\eqref{sh.ex2.} it follows that $V$ is also
indecomposable of degree zero. From the classification,
due to Atiyah, of indecomposable vector bundles of
rank two and degree
zero over an elliptic curve we know that $V$ fits in a
short exact sequence of vector bundles
\begin{equation}\label{sh.ex3.}
0\, \longrightarrow\, L \, \longrightarrow\, V\,
\longrightarrow\, L \, \longrightarrow\, 0\, ,
\end{equation}
where $L$ is a line bundle over $X$ of degree zero
(see \cite[page 432, Theorem 5(ii)]{At}). Let
\begin{equation}\label{sh.ex4.}
0\, \longrightarrow\, F^*_X L \, \longrightarrow\, F^*_XV
\, =\, E\,\longrightarrow\,  F^*_X L \, \longrightarrow\, 0
\end{equation}
be the pull back of the exact sequence in Eq.
\eqref{sh.ex3.} be the Frobenius morphism.

The subbundle ${\mathcal O}_X$ of $E$ in Eq.
\eqref{sh.ex.} is the unique line subbundle of degree
zero. Indeed, if $L'$ is a line subbundle of $E$
of degree zero, then consider the composition
\begin{equation}\label{sh.ex5.}
L' \, \hookrightarrow\, E \, \stackrel{f}{\longrightarrow}
\, {\mathcal O}_X\, ,
\end{equation}
where $f$ is the projection in Eq. \eqref{sh.ex.}. Since
both $L'$ and ${\mathcal O}_X$ are of degree zero, any
nonzero ${\mathcal O}_X$--linear homomorphism between
them must be an isomorphism. Therefore, if $L'$ is
different from the subbundle ${\mathcal O}_X$ of $E$ in Eq.
\eqref{sh.ex.}, then the composition in Eq. \eqref{sh.ex5.}
yields a splitting of the short exact sequence in
Eq. \eqref{sh.ex.}. But this short exact sequence does not
split because $\omega$ in Eq. \eqref{om.} is nonzero.
Hence the subbundle ${\mathcal O}_X$ of $E$ in Eq.
\eqref{sh.ex.} is the unique line subbundle of degree
zero.

Therefore, the subbundle $F^*_X L\, \subset\, E$ in
Eq. \eqref{sh.ex4.} coincides with the
subbundle ${\mathcal O}_X$ of $E$ in Eq. \eqref{sh.ex.}.
This implies that $\omega$ in Eq. \eqref{om.} coincides
with $F^*_X \omega'$, where
$$
\omega'\, \in\, H^1(X,\, {\mathcal O}_X)
$$
is the extension class for the short exact sequence
in Eq. \eqref{sh.ex3.}. But this is impossible because
$X$ is supersingular. Therefore, we conclude that the
vector bundle $E$ is Eq. \eqref{sh.ex.} is not flat.

\subsection{Product of varieties}

Let $X$ and $Y$ be irreducible smooth proper varieties defined
over $k$. Fix a $k$--rational point $x_0$ (respectively,
$y_0$) of $X$ (respectively, $Y$). Recall stratified
group--schemes constructed in Definition \ref{def3}.

\begin{theorem}\label{thm1}
There is a natural isomorphism of the stratified
group--scheme ${\mathcal S}(X\times Y,(x_0,y_0))$
with ${\mathcal S}(X,x_0)\times {\mathcal S}(Y,y_0)$.
\end{theorem}

\begin{proof}
Let
$$
q_X\, :\, X\, \longrightarrow\, X\times Y
$$
be the morphism defined by $x\, \longmapsto\, (x\, ,y_0)$.
Similarly, define
\begin{equation}\label{qY}
q_Y\, :\, Y\, \longrightarrow\, X\times Y
\end{equation}
to be the map that sends any $y$ to $(x_0\, ,y)$.
Consider the corresponding homomorphisms of group--schemes
$$
q^*_X \, :\, {\mathcal S}(X\times Y,(x_0,y_0))\,\longrightarrow\,
{\mathcal S}(X,x_0)
$$
and
$$
q^*_Y \, :\, {\mathcal S}(X\times Y,(x_0,y_0))\,\longrightarrow\,
{\mathcal S}(Y,y_0)
$$
(see Eq. \eqref{e-v}). We will show that the homomorphism
\begin{equation}\label{eta}
\eta\, :=\, q^*_X\times q^*_Y\, :\,
{\mathcal S}(X\times Y,(x_0,y_0))\,\longrightarrow\, 
{\mathcal S}(X,x_0)\times{\mathcal S}(Y,y_0)
\end{equation}
is an isomorphism.

Let
\begin{equation}\label{pX}
p_X\, :\, X\times Y \,\longrightarrow\, X
\end{equation}
and
\begin{equation}\label{pY}
p_Y\, :\, X\times Y \,\longrightarrow\, Y
\end{equation}
be the natural projections. We note that the
composition homomorphism
$$
q^*_X\circ p^*_X\, :\, {\mathcal S}(X,x_0)\,\longrightarrow\,
{\mathcal S}(X,x_0)
$$
is the identity map. Similarly, $q^*_Y\circ p^*_Y$ is also
the identity map. These imply in particular that the homomorphism
$\eta$ in Eq. \eqref{eta} is surjective.

To complete the proof of the theorem we need to show that
$\eta$ is a closed embedding. We will use the criterion in
\cite[page 139, Proposition 2.21(b)]{DM} for a homomorphism
of group--schemes to be a closed embedding.

Let $(V\, ,\nabla)$ be any stratified vector bundle over
$X\times Y$. To show that $\eta$ in Eq. \eqref{eta}
is a closed embedding it suffices to prove that there
are stratified vector bundles
$(E\, ,\nabla^E)$ and $(F\, ,\nabla^F)$ over $X$ and
$Y$ respectively, such that $(V\, ,\nabla)$ is a
quotient of the stratified vector bundle $(p^*_X E\, ,
p^*_X \nabla^E)\bigotimes (p^*_Y F\, ,p^*_Y \nabla^F)$
over $X\times Y$, where $p_X$ and $p_Y$ are the maps
in Eq. \eqref{pX} and Eq. \eqref{pY} respectively.
(See \cite[page 139, Proposition 2.21(b)]{DM}.)

To prove the existence of $(E\, ,\nabla^E)$ and
$(F\, ,\nabla^F)$ satisfying the above condition, consider
the stratified vector bundle
\begin{equation}\label{e4}
(\widehat{V}\, , \widehat{\nabla})\, :=\,
(p^*_Yq^*_YV\, , p^*_Yq^*_Y\nabla)
\end{equation}
on $X\times Y$, where $q_Y$ is the map in Eq. \eqref{qY}
and $p_Y$ is the projection in Eq. \eqref{pY}. Let
\begin{equation}\label{W}
W\, :=\, p_{X*} (V\otimes \widehat{V}^*)
\end{equation}
be the direct image on $X$, where $p_X$ is the projection in
Eq. \eqref{pX}. The stratification $\widehat{\nabla}$
on $\widehat{V}$ (see Eq. \eqref{e4}) and the
stratification $\nabla$ on $V$ together define a
stratification on $V\bigotimes \widehat{V}^*$.
Therefore, we have an induced stratification
\begin{equation}\label{e5}
\widetilde{\nabla}\, :\, 
{\mathcal D}(X)\, \longrightarrow\, {\mathcal E}nd_k(W)
\end{equation}
on the direct image $W$ constructed in Eq. \eqref{W}
(see Eq. \eqref{d.c} for stratification on
a direct image).

Let
\begin{equation}\label{e6}
(U\, ,\nabla^U)\, :=\, (p^*_X W\, ,p^*_X\widetilde{\nabla})
\otimes (\widehat{V}\, , \widehat{\nabla})
\end{equation}
be the tensor product of stratified vector bundles on $X\times Y$,
where $(\widehat{V}\, , \widehat{\nabla})$ is constructed
in Eq. \eqref{e4}, and $\widetilde{\nabla}$ is 
constructed
in Eq. \eqref{e5}. From the construction of
$(U\, ,\nabla^U)$ in Eq. \eqref{e6} it follows that
we have a homomorphism of stratified vector bundles
\begin{equation}\label{e7}
\nu\, :\, (U\, ,\nabla^U)\, \longrightarrow\, 
(V\, ,\nabla)\, .
\end{equation}
To explain this with more detail, we
note that for any two ${\mathcal O}_{X\times Y}$--coherent
sheaves $A$ and $B$ on $X\times Y$, we have a
natural homomorphism
\begin{equation}\label{n.h.}
p^*_Xp_{X*}(B\otimes A^*)\otimes A\, \longrightarrow\, B
\end{equation}
of ${\mathcal O}_{X\times Y}$--coherent sheaves which is
constructed using the obvious homomorphism
$$
p^*_Xp_{X*}(B\otimes A^*)\,\longrightarrow\,B\otimes A^*\, .
$$
Let 
$$
U\, =\, p^*_Xp_{X*} (V\otimes (p^*_Yq^*_YV)^*)\otimes
p^*_Yq^*_YV \, \longrightarrow\, V
$$
be the homomorphism of
${\mathcal O}_{X\times Y}$--coherent sheaves
in Eq. \eqref{n.h.} obtained by substituting
$p^*_Yq^*_YV$ for $A$ and $V$ for $B$. This
homomorphism intertwines $\nabla^U$ (defined in
Eq. \eqref{e6}) and $\nabla$, giving the homomorphism
$\nu$ in Eq. \eqref{e7}.

Since $Y$ is proper, using
\cite[page 9, Proposition 2.4]{Gi} it follows immediately
that the homomorphism $\nu$ in Eq. \eqref{e7} is surjective.
Therefore, if we set
$$
(E\, ,\nabla^E)\, :=\, (W\, ,\widetilde{\nabla})
$$
where $W$ and $\widetilde{\nabla}$ are constructed in
Eq. \eqref{W} and Eq. \eqref{e5} respectively, and also set
$$
(F\, ,\nabla^F)\, :=\, (q^*_YV\, , q^*_Y\nabla)
$$
(see Eq. \eqref{e4}), then
$$
(U\, ,\nabla^U)\, =\, (p^*_X E\, ,
p^*_X \nabla^E)\otimes (p^*_Y F\, ,p^*_Y \nabla^F)\, .
$$
Since the homomorphism $\nu$ in Eq. \eqref{e7} is surjective,
we now conclude that the stratified vector bundle $(V\, ,
\nabla)$ is a quotient of the tensor product $(p^*_X E\, ,
p^*_X \nabla^E)\bigotimes (p^*_Y F\, ,p^*_Y \nabla^F)$.
We noted earlier that this implies that the homomorphism
$\eta$ in Eq. \eqref{eta} is a closed embedding.

Therefore, $\eta$ is an isomorphism. This completes the
proof of the theorem.
\end{proof}

\section{Automorphisms of the group--scheme}

As in Section \ref{sec2.1}, we take $X$ to be an
irreducible smooth variety defined over $k$.
Since the field $k$ is algebraically closed,
in particular $k$ is perfect, the arithmetic
Frobenius homomorphism $k\, \longrightarrow\, k$
is an isomorphism. Hence we can
identify all $X_n$, $n\, \geq\, 1$, with $X$ itself.
Consequently, the Frobenius map $F_X$ becomes a
self--map of $X$.

Given a flat vector bundle $\{E_n\, ,\sigma_n\}_{n\geq 0}$
on $X$ (see Definition \ref{def2}), we may construct a
new flat vector bundle on $X$ as follows. For each
$n\,\geq\, 0$, set
$$
E'_n\, :=\, E_{n+1}\, .
$$
Set
$$
\sigma'_n \, :=\, \sigma_{n+1}\, :\,
F^*_XE'_{n+1}\, =\, F^*_XE_{n+2}\, \longrightarrow\,
E_{n+1}\, =\, E'_n
$$
to be the isomorphism.
Now $\{E'_n\, ,\sigma'_n\}_{n\geq 0}$ is a new flat
vector bundle on $X$.

Using \cite[page 4, Theorem 1.3]{Gi}, which identifies
flat vector bundles on $X$ with stratified vector bundles
on $X$, the above construction produces a functor from the
category ${\mathcal C}(X)$
of stratified vector bundles on $X$ to itself.
In fact, this construction yields a functor from the
neutral Tannakian category $({\mathcal C}(X)\, , T_{x_0})$
over $k$ (see Definition \ref{def3}) to itself. Therefore,
we get an endomorphism of the stratified group--scheme
\begin{equation}\label{F}
\Phi\, :\, {\mathcal S}(X,x_0)\, \longrightarrow\,
{\mathcal S}(X,x_0)
\end{equation}
(see Definition \ref{def3} for ${\mathcal S}(X,x_0)$).

\begin{lemma}\label{lem3}
The homomorphism $\Phi$ in Eq. \eqref{F} is an isomorphism.
\end{lemma}

\begin{proof}
Let $\{E_n\, ,\sigma_n\}_{n\geq 0}$ be a flat vector bundle
on $X$. We will construct another flat vector bundle from it.
For $n\, \geq\, 0$, set
$$
F_n\, :=\, F^*_X E_n\, .
$$
So for $n\, \geq\, 1$, using the isomorphism $\sigma_{n-1}$,
the vector bundle $F_n$ gets identified with $E_{n-1}$.
For $n\, \geq\, 0$, set
$$
\tau_n\, :=\, F^*_X\sigma_n\, :\, F^*_X F_{n+1}\, =\,
F^*_X(F^*_X E_{n+1}) \, \longrightarrow\, F^*_X E_n\, =\, F_n\, .
$$
It is easy to check that $\{F_n\, ,\tau_n\}_{n\geq 0}$ is a flat 
vector bundle on $X$.

Using the identification between stratified vector bundles on $X$
and flat vector bundles on $X$, the above construction produces
a functor from the neutral Tannakian category
$({\mathcal C}(X)\, , T_{x_0})$ over $k$ in Definition \ref{def3}
to itself. Consequently, we get a
homomorphism of group--schemes
\begin{equation}\label{phi-i}
\Psi\, :\, {\mathcal S}(X,x_0)\, \longrightarrow\,
{\mathcal S}(X,x_0)\, .
\end{equation}
It is now straight--forward to check that
\begin{equation}\label{phi-i2}
\Psi\circ \Phi\, =\, \Phi\circ\Psi\, =\, \text{Id}_{{\mathcal 
S}(X,x_0)}\, .
\end{equation}
This completes the proof of the lemma.
\end{proof}

{}From the construction of the homomorphism $\Psi$ in Eq.
\eqref{phi-i} it follows that $\Psi$
coincides with the Frobenius morphism
$$
F_{{\mathcal S}(X,x_0)}\, :\,
{\mathcal S}(X,x_0)\, \longrightarrow\,
{\mathcal S}(X,x_0)
$$
of the group--scheme ${\mathcal S}(X,x_0)$ (see
\cite[page 146--148]{Ja} for the Frobenius morphism of a
group--scheme). We have seen in Eq. \eqref{phi-i2} that
$\Psi$ is an isomorphism. Hence the Frobenius
morphism $F_{{\mathcal S}(X,x_0)}$ is an isomorphism.

We put down this observation as the following corollary.

\begin{corollary}\label{cor.fi}
The Frobenius morphism $F_{{\mathcal S}(X,x_0)}$
of the group--scheme ${\mathcal S}(X,x_0)$ is an
isomorphism.
\end{corollary}

\section{Flat principal bundles}

\subsection{A tautological principal bundle}\label{se5.1}

Consider the group--scheme ${\mathcal S}(X,x_0)$ in
Definition \ref{def3}. We will show that there is
a tautological principal ${\mathcal S}(X,x_0)$--bundle
over $X$. 

Let $\text{Vect}(X)$ denote the category of vector bundles
over $X$. We have a functor
\begin{equation}\label{B}
{\mathcal B}\, :\, {\mathcal C}(X)\,\longrightarrow\,
\text{Vect}(X)\, ,
\end{equation}
where ${\mathcal C}(X)$ as before is the category of
stratified vector bundles on $X$,
that sends any stratified vector bundle $(E\, ,\nabla)$ over
$X$ to the vector bundle $E$. This functor defines a principal
${\mathcal S}(X,x_0)$--bundle over $X$; see
\cite[Lemma 2.3, Proposition 2.4]{No} (reproduced as Theorem 2.3
in \cite[page 6]{BPS}), \cite[page 149, Theorem 3.2]{DM}.
Let
\begin{equation}\label{e8}
E_{{\mathcal S}(X,x_0)}\, \longrightarrow\, X
\end{equation}
be the principal ${\mathcal S}(X,x_0)$--bundle given by
the functor $\mathcal B$ in Eq. \eqref{B}.

Therefore, finite dimensional representations
of the affine group--scheme
${\mathcal S}(X,x_0)$ are stratified vector bundles on $X$,
and furthermore, the vector bundle over $X$ associated to the
principal ${\mathcal S}(X,x_0)$--bundle
$E_{{\mathcal S}(X,x_0)}$ in Eq. \eqref{e8} for a
representation $(E\, ,\nabla)$ of ${\mathcal S}(X,x_0)$
is the vector bundle $E$ itself.

\subsection{Flat principal bundles}

Let $\mathcal G$ be an affine group--scheme
defined over the field $k$.
Imitating the definition of a flat vector bundle we define:

\begin{definition}\label{def4}
{\rm A} flat principal $\mathcal G$--bundle {\rm over $X$ is a
sequence of pairs $\{E_n\, ,\sigma_n\}_{n\geq 0}$, where $E_n$
is a principal $\mathcal G$--bundle over $X_n$ and}
$$
\sigma_n\, :\, F^*_X E_{n+1}\, \longrightarrow\, E_{n}
$$
{\rm is an isomorphism of principal $\mathcal G$--bundles.}

{\rm A flat principal $\mathcal G$--bundle $\{E_n\, ,
\sigma_n\}_{n\geq 0}$ will also be called a} flat structure
{\rm on the principal $\mathcal G$--bundle $E_0$.}
\end{definition}

\begin{lemma}\label{lem4}
The principal ${\mathcal S}(X,x_0)$--bundle
$E_{{\mathcal S}(X,x_0)}$ in Eq. \eqref{e8} has a natural flat
structure.
\end{lemma}

\begin{proof}
For any positive integer $n$, let
\begin{equation}\label{phi-n}
\Phi^n \, :=\, \overbrace{\Phi\circ\cdots\circ 
\Phi}^{n\mbox{-}\rm{times}}
\, :\, {\mathcal S}(X,x_0)\, \longrightarrow\,
{\mathcal S}(X,x_0)
\end{equation}
be the $n$--fold iteration of the homomorphism
$\Phi$ constructed in Eq. \eqref{F}, and set $\Phi^0$
to be the identity map of ${\mathcal S}(X,x_0)$.
For $n\, \geq\, 0$, let
\begin{equation}\label{e9}
E^n_{{\mathcal S}(X,x_0)}\, :=\,
E_{{\mathcal S}(X,x_0)}(\Phi^n) \, \longrightarrow\, X
\end{equation}
be the principal ${\mathcal S}(X,x_0)$--bundle over $X$
obtained by extending the structure group of the
principal ${\mathcal S}(X,x_0)$--bundle
$E_{{\mathcal S}(X,x_0)}$ in Eq. \eqref{e8} using the
homomorphism $\Phi^n$ in Eq. \eqref{phi-n}.

Let $F_{{\mathcal S}(X,x_0)}$ be any principal ${\mathcal
S}(X,x_0)$--bundle over $X$. Let
$F_{{\mathcal S}(X,x_0)}(\Psi)$ be the principal ${\mathcal 
S}(X,x_0)$--bundle over $X$ obtained by extending the structure
group of $F_{{\mathcal S}(X,x_0)}$ using the homomorphism
$\Psi$ in Eq. \eqref{phi-i}. From the construction of $\Psi$
it follows that $F_{{\mathcal S}(X,x_0)}(\Psi)$ is
canonically identified the pull back
$F^*_X F_{{\mathcal S}(X,x_0)}$.
Using this identification for the principal bundle
$$
F_{{\mathcal S}(X,x_0)}\, =\, E^{n+1}_{{\mathcal S}(X,x_0)}
$$
defined in Eq. \eqref{e9}, together with the fact that
$$
\Psi\circ \Phi^{n+1}\, =\, \Phi^n\, ,
$$
which follows from the observation in
the proof of Lemma \ref{lem3} that
$$
\Psi\circ \Phi\, =\, \text{Id}_{{\mathcal S}(X,x_0)}\, ,
$$
we get an isomorphism of the principal ${\mathcal
S}(X,x_0)$--bundle $F^*_X E^{n+1}_{{\mathcal S}(X,x_0)}$
with $E^n_{{\mathcal S}(X,x_0)}$.

Let
$$
\sigma_n\, :\, F^*_X E^{n+1}_{{\mathcal S}(X,x_0)}\,
\longrightarrow\, E^n_{{\mathcal S}(X,x_0)}
$$
be the isomorphism of principal ${\mathcal
S}(X,x_0)$--bundles obtained above. Therefore,
\begin{equation}\label{fi2.}
\{E^n_{{\mathcal S}(X,x_0)}\, ,\sigma_n\}_{n\geq 0}
\end{equation}
is a flat principal ${\mathcal S}(X,x_0)$--bundle over $X$.
This completes the proof of the lemma.
\end{proof}

\subsection{Universality of the tautological principal bundle}

Let $\mathcal G$ be any affine group--scheme defined over
$k$. Let $\{E_n\, ,\sigma_n\}_{n\geq 0}$ be a flat principal
$\mathcal G$--bundle
over an irreducible smooth variety $X$ defined over $k$;
see Definition \ref{def4} for flat $\mathcal G$--bundle. Let

$$
\text{Ad}(E_0)\, =\, E_0({\mathcal G})
$$
be the adjoint bundle of the principal $\mathcal G$--bundle
$E_0$. We recall that the adjoint bundle $\text{Ad}(E_0)$
is the fiber bundle over $X$ associated to $E_0$ for the
adjoint action of $\mathcal G$ on itself. So $\text{Ad}(E_0)$
is a group--scheme over $X$.

Let $\text{Ad}(E_0)_{x_0}$ be the fiber of $\text{Ad}(E_0)$ over
the base point $x_0$. So the group--scheme $\text{Ad}(E_0)_{x_0}$
is isomorphic to $\mathcal G$.

Fix a $k$--rational point
\begin{equation}\label{base2}
\widehat{x}_0\, \in\, (E_0)_{x_0}
\end{equation}
in the fiber of $E_0$ of $x_0$. Since the
pull back of any principal bundle to the total space of the
principal bundle is canonically trivialized,
using the base point $\widehat{x}_0$ in Eq. \eqref{base2}
we get an isomorphism of group--schemes
\begin{equation}\label{base3}
\rho_0\, :\, \text{Ad}(E_0)_{x_0}\, \longrightarrow
\, {\mathcal G}
\end{equation}
(see also \cite[Section~3]{BPS}).

\begin{proposition}\label{prop2}
There is homomorphism of group schemes
$$
\rho\, :\, {\mathcal S}(X,x_0)\, \longrightarrow\,
{\rm Ad}(E_0)_{x_0}
$$
canonically associated to the flat principal
$\mathcal G$--bundle $\{E_n\, ,\sigma_n\}_{n\geq 0}$.

The flat principal $\mathcal G$--bundle $\{E_n\, ,
\sigma_n\}_{n\geq 0}$ coincides with the one obtained by
extending the structure group of the flat principal ${\mathcal 
S}(X,x_0)$--bundle $\{E^n_{{\mathcal S}(X,x_0)}\, ,
\sigma_n\}_{n\geq 0}$ (see Lemma \ref{lem4}) using the
homomorphism $\rho_0\circ \rho$, where $\rho_0$ is
constructed in Eq. \eqref{base3}.
\end{proposition}

\begin{proof}
We will first recall a Tannakian description of the
group--scheme $\text{Ad}(E_0)_{x_0}$.

Let $\text{Rep}({\mathcal G})$ denote the category of all
finite dimensional representations of the group--scheme
${\mathcal G}$. It is a rigid abelian
$k$--linear tensor category. We will construct a fiber
functor on $\text{Rep}({\mathcal G})$.
For any finite dimensional left $\mathcal G$--module $V$, let
$E_0(V)$ denote the vector bundle over $X$ associated to
the principal ${\mathcal G}$--bundle $E_0$ for the
${\mathcal G}$--module $V$. Let $E_0(V)_{x_0}$ be the
fiber of $E_0(V)$ over the base point $x_0$. Now we
have a fiber functor
\begin{equation}\label{fi.}
T_{E_0}\, :\, \text{Rep}({\mathcal G})\, \longrightarrow\,
\text{Vect}(k)
\end{equation}
that sends any $V$ to $E_0(V)_{x_0}$. So the pair
$(\text{Rep}({\mathcal G})\, , T_{E_0})$ defines a neutral
Tannakian category over $k$. The corresponding affine
group--scheme over $k$ (see \cite[page 130, Theorem
2.11]{DM}) is identified with $\text{Ad}(E_0)_{x_0}$.

Take any $V\, \in\,
\text{Rep}({\mathcal G})$. For each $n\, \geq\, 0$,
let $E_n(V)$ denote the vector bundle over $X$ associated to
the principal ${\mathcal G}$--bundle $E_n$ for the ${\mathcal
G}$--module $V$. The isomorphism of principal
${\mathcal G}$--bundles
$$
\sigma_n\, :\, F^*_X E_{n+1}\, \longrightarrow\, E_{n}
$$
induces an isomorphism of associated vector bundles
$$
\widehat{\sigma}_n\, :\, F^*_X E_{n+1}(V)\, \longrightarrow\,
E_{n}(V)\, .
$$

It is straight--forward to check that
\begin{equation}\label{f.b.}
\{E_n(V)\, , \widehat{\sigma}_n\}_{n\geq 0}
\end{equation}
is a flat vector bundle over $X$.
Therefore, $\{E_n(V)\, ,\widehat{\sigma}_n\}_{n\geq 0}$ gives a 
stratified vector bundle over $X$. Let
$$ 
(V'\, ,\nabla')\, \in\, {\mathcal C}(X)
$$
be the stratified vector bundles given by the flat vector bundle 
$\{E_n(V)\, ,\widehat{\sigma}_n\}_{n\geq 
0}$. Consequently, we have a functor
\begin{equation}\label{fun.}
\text{Rep}({\mathcal G})\, 
\longrightarrow\, {\mathcal C}(X)
\end{equation}
that sends any $V$ to $(V'\, 
,\nabla')$. Now, comparing the two fiber functors $T_{x_0}$
and $T_{E_0}$, defined in Eq. \eqref{e1} and Eq. \eqref{fi.}
respectively, we see that the functor in Eq. \eqref{fun.}
actually produces to a functor from the neutral Tannakian
category $(\text{Rep}({\mathcal G})\, ,T_{E_0})$ over $k$
to the neutral Tannakian category $({\mathcal C}(X)\, ,
T_{x_0})$. In view of the above Tannakian description of the
group--scheme $\text{Ad}(E_0)_{x_0}$, this functor between
neutral Tannakian categories over $k$ produces a homomorphism
\begin{equation}\label{rho.f}
\rho\, :\, {\mathcal S}(X,x_0)\, \longrightarrow\,
\text{Ad}(E_0)_{x_0}
\end{equation}
of group--schemes over $k$.

To prove the second part of the proposition, for each
$n\, \geq\, 0$, let
$$
F^n_{\mathcal G}\, :=\, E^n_{{\mathcal S}(X,x_0)}({\mathcal G})
$$
be the principal $\mathcal G$--bundle over $X$ obtained by
extending the principal ${\mathcal S}(X,x_0)$--bundle
$E^n_{{\mathcal S}(X,x_0)}$ (see Eq. \eqref{fi2.})
using the homomorphism $\rho_0\circ \rho$, where $\rho$ is
the homomorphism constructed in Eq. \eqref{rho.f}, and
$\rho_0$ is the homomorphism in Eq. \eqref{base3}. Let
\begin{equation}\label{tau.n}
\tau_n\, :\, F^*_XF^{n+1}_{\mathcal G}\, \longrightarrow\,
F^n_{\mathcal G}
\end{equation}
be the isomorphism of principal $\mathcal G$--bundles
induced by the isomorphism $\sigma_n$ in Eq. \eqref{fi2.}.
Therefore, $\{F^n_{\mathcal G}\, ,\tau_n\}_{n\geq 0}$
is a flat principal $\mathcal G$--bundle over $X$.

Take any $V\, \in\, \text{Rep}({\mathcal G})$. For
each $n\, \geq\, 0$, let
$$
F^n_V\, :=\, F^n_{\mathcal G}(V)
$$
be the vector bundle over $X$ associated to the principal
$\mathcal G$--bundle $F^n_{\mathcal G}$ for the $G$--module
$V$. Let
$$
\widehat{\tau}_n\, :\, F^*_X F^{n+1}_V \, \longrightarrow\,
F^n_V
$$
be the isomorphism of vector bundles induced by the isomorphism
${\tau}_n$ in Eq. \eqref{tau.n}. From the construction of
the homomorphism $\rho$ in Eq. \eqref{rho.f} it follows that
the flat vector bundle
$$
\{F^n_V\, ,\widehat{\tau}_n\}
$$
is naturally identified with the flat vector bundle $\{E_n(V)
\, , \widehat{\sigma}_n\}_{n\geq 0}$ in Eq. \eqref{f.b.}. From
this it follows that  the flat
principal $\mathcal G$--bundle $\{E_n\, ,
\sigma_n\}_{n\geq 0}$ coincides with
$\{F^n_{\mathcal G}\, ,\tau_n\}_{n\geq 0}$ (see
\cite[page 149, Theorem 3.2]{DM}
\cite[Lemma 2.3, Proposition 2.4]{No}).
This completes the proof of the proposition.
\end{proof}

%%%%%%%%%%%%%%%%%%%%%%%%%%%%%%%%%%%%%%%%%%%%%%%%%%%%%%%%%%%%%%

\end{document}